\documentclass[12pt]{article}
\usepackage{graphicx}
\usepackage{color}
\catcode`\@=11 \@addtoreset{equation}{section}

\catcode`\@=12

\newtheorem{Theorem}{Theorem}[section]
\newtheorem{Proposition}{Proposition}[section]
\newtheorem{Lemma}{Lemma}[section]
\newtheorem{Corollary}{Corollary}[section]
\newtheorem{Remark}{Remark}[section]
\newtheorem{Definition}{Definition}[section]

\newcommand{\bTheorem}[1]{
\begin{Theorem} \label{T#1} }
\newcommand{\eT}{\end{Theorem}}

\newcommand{\bProposition}[1]{
\begin{Proposition} \label{P#1}}
\newcommand{\eP}{\end{Proposition}}

\newcommand{\bLemma}[1]{
\begin{Lemma} \label{L#1} }
\newcommand{\eL}{\end{Lemma}}

\newcommand{\bCorollary}[1]{
\begin{Corollary} \label{C#1} }
\newcommand{\eC}{\end{Corollary}}

\newcommand{\bDefinition}[1]{
\begin{Definition} \label{D#1} }
\newcommand{\eD}{\end{Definition}}

\newcommand{\bRemark}[1]{
\begin{Remark} \label{R#1} }
\newcommand{\eR}{\end{Remark}}

\newcommand{\bFormula}[1]{
\begin{equation} \label{#1}}
\newcommand{\eF}{\end{equation}}

\newcommand{\Ov}[1]{\overline{#1}}
\newcommand{\res}[1] { \left[ #1 \right]_{\rm res} }
\newcommand{\ess}[1] { \left[ #1 \right]_{\rm ess} }

\newcommand{\DC}{C^\infty_c}

\newcommand{\vr}{\varrho}
\newcommand{\vre}{\vr_\ep}

\newcommand{\vt}{\vartheta}
\newcommand{\vu}{\vc{u}}
\newcommand{\vc}[1]{{\bf #1}}

\newcommand{\Div}{{\rm div}_x}
\newcommand{\Grad}{\nabla_x}

\newcommand{\tn}[1]{\mbox {\F #1}}
\newcommand{\dx}{{\rm d} {x}}
\newcommand{\dt}{{\rm d} t }

\newcommand{\intO}[1]{\int_{\Omega} #1 \ \dx}

\newcommand{\aleq}{\stackrel{<}{\sim}}
\newcommand{\ageq}{\stackrel{>}{\sim}}

\newcommand{\ep}{\varepsilon}

\font\F=msbm10 scaled 1000

\newcommand{\vuE}{\vu_E}
\newcommand{\vrE}{\vr_E}
\newcommand{\vtE}{\vt_E}

\definecolor{Cgrey}{rgb}{0.85,0.85,0.85}
\definecolor{Cblue}{rgb}{0.50,0.85,0.85}
\definecolor{Cred}{rgb}{1,0,0}
\definecolor{fancy}{rgb}{0.10,0.85,0.10}

\newcommand\Cbox[2]{%
    \newbox\contentbox%
    \newbox\bkgdbox%
    \setbox\contentbox\hbox to \hsize{%
        \vtop{
            \kern\columnsep
            \hbox to \hsize{%
                \kern\columnsep%
                \advance\hsize by -2\columnsep%
                \setlength{\textwidth}{\hsize}%
                \vbox{
                    \parskip=\baselineskip
                    \parindent=0bp
                    #2
                }%
                \kern\columnsep%
            }%
            \kern\columnsep%
        }%
    }%
    \setbox\bkgdbox\vbox{
        \color{#1}
        \hrule width  \wd\contentbox %
               height \ht\contentbox %
               depth  \dp\contentbox
        \color{black}
    }%
    \wd\bkgdbox=0bp%
    \vbox{\hbox to \hsize{\box\bkgdbox\box\contentbox}}%
    \vskip\baselineskip%
}

\date{}

\begin{document}


\title{Vanishing dissipation limit for the Navier-Stokes-Fourier system}

\author{Eduard Feireisl \thanks{The research of E.F. leading to these results has received funding from the European Research Council under the European Union's Seventh Framework
Programme (FP7/2007-2013)/ ERC Grant Agreement 320078. The Institute of Mathematics of the Academy of Sciences of the Czech
Republic is supported by RVO:67985840.
 } }

\maketitle

\bigskip

\centerline{Institute of Mathematics of the Czech Academy of Sciences}

\centerline{\v Zitn\' a 25, 115 67 Praha 1, Czech Republic}






\bigskip





\begin{abstract}

We consider the motion of a compressible, viscous, and heat conducting fluid in the regime of small viscosity and heat conductivity.
It is shown that weak solutions of the associated Navier-Stokes-Fourier system converge to
a (strong) solution of the Euler system on its life span. The problem is studied in a bounded domain $\Omega \subset R^3$, on the boundary of which the
velocity field satisfies the complete slip boundary conditions.

\end{abstract}



\section{Introduction}
\label{i}

A formal limit of vanishing viscosity and heat conductivity in the fluid models based on the principles of continuum mechanics gives rise to the
\emph{Euler system}
\bFormula{i1}
\partial_t \vr + \Div (\vr \vu) = 0,
\eF
\bFormula{i2}
\partial_t (\vr \vu) + \Div (\vr \vu \otimes \vu) + \Grad p_M(\vr, \vt) = 0,
\eF
\bFormula{i3}
\partial_t \left( \frac{1}{2} \vr |\vu|^2 + \vr e_M(\vr, \vt) \right) + \Div \left[ \left( \frac{1}{2} \vr |\vu|^2 + \vr e_M(\vr, \vt) \right) \vu +
p_M (\vr, \vt) \vu \right] = 0,
\eF
describing the time evolution of the basic macroscopic quantities: the mass density $\vr = \vr(t,x)$, the absolute temperature $\vt = \vt(t,x)$, and the velocity field $\vu = \vu(t,x)$. The
symbols $p_M$, $e_M$ denote the (molecular) pressure and the associated specific internal energy, respectively. Solutions of system (\ref{i1}--\ref{i3})
are known to develop singularities in a finite time lap even if the initial state
\bFormula{i4}
\vr(0, \cdot) = \vr_0, \ \vt(0, \cdot) = \vt_0, \ \vu(0, \cdot) = \vu_0
\eF
is regular, represented by smooth functions. Accordingly, solutions of (\ref{i1}--\ref{i2}) are usually understood in the weak sense, where all derivatives are interpreted as mathematical distributions.

Unfortunately, the class of weak solutions is too large to secure uniqueness and/or
continuous dependence of solutions on the data. To remedy this problem, several additional \emph{admissibility criteria} have been proposed, among which the
entropy inequality
\bFormula{i5}
\partial_t (\vr s_M(\vr, \vt)) + \Div (\vr s_M(\vr, \vt) \vu) \geq 0,
\eF
where $s_M = s_M(\vr, \vt)$ is the specific entropy related to $e_M$ and $p_M$ through Gibbs' equation
\bFormula{i6}
\vt D s_M(\vr, \vt) = D e_M(\vr, \vt) + p_M(\vr, \vt) D \left( \frac{1}{\vr} \right).
\eF
Although the entropy and similar admissibility criteria based on the Second law of thermodynamics have been partially successful when applied to problems
in the 1-D geometry, see Bressan \cite{BRESSAN}, Dafermos \cite{D4}, they failed in identifying the relevant solution in the natural 3-D setting, see
DeLellis and Sz\' ekelyhidi \cite{DelSze3}, Chiodaroli, DeLellis, and Kreml \cite{ChiDelKre}, Chiodaroli and Kreml \cite{ChiKre}

In the light of the above arguments and in accordance with the general approach advocated by Bardos et al. \cite{BarTi}, the physically relevant
solutions to the Euler system (\ref{i1}--\ref{i3}) should be identified as asymptotic limits of solutions to more complex problems - primitive systems - describing the evolution of ``real'' fluids. As an example of such primitive problem, we consider the \emph{Navier-Stokes-Fourier system} in the form
\bFormula{e1}
\partial_t \vr + \Div (\vr \vu) = 0,
\eF
\bFormula{e2}
\partial_t (\vr \vu) + \Div (\vr \vu \otimes \vu) + \Grad p(\vr, \vt) = \Div \tn{S}(\vt, \Grad \vu) - \lambda \vu
\eF
\bFormula{e3}
\partial_t (\vr s(\vr, \vt)) + \Div (\vr s(\vr, \vt) \vu) + \Div \left( \frac{\vc{q}}{\vt} \right) =
\sigma, \ \sigma = \frac{1}{\vt} \left( \tn{S}(\vt, \Grad \vu) - \frac{\vc{q} \cdot \Grad \vt}{\vt} \right),
\eF
where $\tn{S}(\vt, \Grad \vu)$ is the viscous stress tensor given by Newton's law
\bFormula{e5}
\tn{S}(\vt, \Grad \vu) = \nu \left[ \mu(\vt) \left( \Grad \vu + \Grad^t \vu - \frac{2}{3} \Div \vu \tn{I} \right) + \eta(\vt) \Div \vu \tn{I} \right],
\ \nu > 0,
\eF
and $\vc{q} = \vc{q}(\vt, \Grad \vt)$ is the heat flux determined by Fourier's law
\bFormula{e6}
\vc{q} = - \omega \kappa(\vt) \Grad \vt,\ \omega > 0.
\eF

The pressure $p = p_M + p_R$ is augmented by a radiation component,
\bFormula{e7}
p(\vr, \vt) = p_M(\vr, \vt) + p_R(\vr, \vt), \ p_R(\vr, \vt) = \frac{a}{3} \vt^4,\ a > 0,
\eF
while the internal energy reads
\bFormula{e8}
e(\vr, \vt) = e_M (\vr, \vt) + e_R(\vr, \vt),\
e_R = a \frac{\vt^4}{\vr}.
\eF
Accordingly, the specific entropy is
\bFormula{e8a}
s(\vr, \vt) = s_M (\vr, \vt) + s_R(\vr, \vt), \ s_R(\vr, \vt) = s_R (\vr, \vt) = \frac{4a}{3 \vr} \vt^3,
\eF
see \cite[Chapter 1]{FEINOV} for the physical background.

Our main goal will be to show that solutions of the Navier-Stokes-Fourier system (\ref{e1}--\ref{e3}) converge to those of the Euler system
(\ref{i1}--\ref{i3}) in the asymptotic regime
\bFormula{e8b}
\nu, \ \omega, \ a , \ \lambda \to 0,
\eF
on condition that the limit solution is smooth. In order to avoid the well-known and still unsurmountable difficulties connected with the presence of
a boundary layer (see e.g. the survey by E \cite{E1}), we impose the complete slip boundary conditions for the velocity
\bFormula{e20}
\vc{u} \cdot \vc{n}|_{\partial \Omega} = 0, \ \left[ \tn{S}(\vt, \Grad \vu) \cdot \vc{n} \right] \times \vc{n}|_{\partial \Omega} = 0,
\eF
accompanied with the no-flux condition
\bFormula{e21}
\vc{q}(\vt, \Grad \vt) \cdot \vc{n} |_{\partial \Omega} = 0,
\eF
where $\vc{n}$ denotes the outer normal vector to $\partial \Omega$.

Besides a vast amount of literature concerning the inviscid limit for the \emph{incompressible} Navier-Stokes system, see e.g. Kato \cite{Kato},
Temam, Wang \cite{TemWan1}, \cite{TemWan2}, Wang, Xin, and Zang
\cite{WaXiZa}, the survey articles by E \cite{E1}, Masmoudi \cite{MAS1}, and the references cited therein, much less seems to be known in the context of compressible fluids. There are results on stability of certain special solutions like shock or rarefaction waves, see Hoff and Liu \cite{HofLiu}, Hoff and Zumbrun \cite{HOZU}, Li and Wang \cite{LiWang} as well as studies of the related linearized problems, Xin and Yanagisawa \cite{XinYan}. Probably the closest result to ours has been recently obtained by Sueur \cite{Sue1}
for the barotropic Navier-Stokes system, related issues were discussed by Wang and Williams \cite{WanWil}.

Our approach is based on the \emph{relative energy} (see Dafermos \cite{Daf4} for a closely related concept of relative entropy) and the associated concept of \emph{dissipative solution}  for the full Navier-Stokes-Fourier system introduced in \cite{FeiNov10}. A similar strategy has been already used in \cite{FeiNov12} in the context of inviscid-incompressible limits, and in \cite{FeiNov10}, were the property of weak-strong uniqueness was established. In comparison with these problems, the purely inviscid, or, more precisely, zero dissipation limit features new
additional difficulties:
\begin{itemize}
\item In the inviscid-incompressible limit studied in \cite{FeiNov12}, the initial data are prepared, in particular, the density and the temperature in the primitive system are {\it a priori} known to be close to their limit (constant) values. Such a piece of information is not available for the vanishing dissipation limit.

\item In contrast with the situation in \cite{FeiNov10}, where weak and strong solutions of the same (viscous) problem are being compared, the uniform bounds based on the presence of
viscosity and heat conductivity are lost in the vanishing dissipation limit.
\end{itemize}

As a consequence of the afore mentioned difficulties, our result is \emph{path dependent} - the rates of convergence to zero of the singular parameters $\nu$, $\omega$, $a$, and $\lambda$ are interrelated in a special way specified in Section \ref{m}. For the same reason, the presence of the
``damping'' term $-\lambda \vu$ in the primitive system is necessary to control the amplitude of certain quantities, notably the velocity, on the hypothetical
vacuum zones created by vanishing density.

The paper is organized as follows. In Section \ref{e}, we summarize the necessary preliminary material including the main hypotheses and the concept of
weak solution for the Navier-Stokes-Fourier system. In Section \ref{m}, we state our main result. The
relative functional together with the associated relative energy inequality adapted to the present setting
are introduced in Section \ref{r}. The available uniform bounds on the family of solutions to the primitive system are collected in Section \ref{E}. The proof of convergence is completed in Section \ref{c}. Possible extensions are discussed in Section \ref{F}.

\section{Preliminaries, weak solutions}
\label{e}

We adopt the concept of weak solution to the Navier-Stokes-Fourier system (\ref{e1}--\ref{e6}), (\ref{e20}), (\ref{e21}) introduced in \cite[Chapter 2]{FEINOV}.
More specifically,
\begin{itemize}
\item
the equations (\ref{e1}), (\ref{e2}), together with the boundary conditions (\ref{e20}), are understood in the sense of distributions;
\item the entropy balance (\ref{e3}), with (\ref{e21}), is relaxed to an \emph{inequality}
\bFormula{e3a}
\partial_t (\vr s(\vr, \vt)) + \Div (\vr s(\vr, \vt) \vu) + \Div \left( \frac{\vc{q}}{\vt} \right) =
\sigma, \ \sigma \geq \frac{1}{\vt} \left( \tn{S}(\vt, \Grad \vu) - \frac{\vc{q} \cdot \Grad \vt}{\vt} \right)
\eF
satisfied in the sense of distribution;
\item
the system is augmented by the total energy balance
\bFormula{e4}
\intO{ \left[ \frac{1}{2} \vr |\vu|^2 + \vr e(\vr, \vt) \right] (\tau, \cdot) } + \lambda \int_0^\tau \intO{ |\vu|^2 } \leq
\intO{ \left[ \frac{1}{2} \vr_0 |\vu_0|^2 + \vr_0 e(\vr_0, \vt_0) \right] }
\eF
for a.a. $\tau \in [0,T]$.
\end{itemize}

\bRemark{e1}
As a matter of fact, the original definition in \cite[Chapter 2]{FEINOV} does not include the damping term $-\lambda \vu$, and, accordingly, stipulates  \emph{equality} rather than inequality in the total energy balance (\ref{e4}). In view of the anticipated lack of compactness of the velocity field on the vacuum, however,
the inequality seems more adequate in the present setting.
\eR

The existence theory developed in \cite[Chapter 3]{FEINOV} requires certain structural restrictions to be imposed on the constitutive relations listed below:

\begin{itemize}
\item The pressure $p$ takes the form (\ref{e7}), with
\bFormula{H1}
p_M(\vr, \vt) = \vt^{5/2} P \left( \frac{\vr}{\vt^{3/2}} \right),
\eF
where $P \in C^1 [0, \infty) \cap C^5(0, \infty)$ satisfies
\bFormula{H2}
P(0) = 0, \ P'(Z) > 0 \ \mbox{for all} \ Z \geq 0,
\eF
\bFormula{H3}
0 < \frac{ \frac{5}{3} P(Z) - P'(Z) Z }{Z} < c \ \mbox{for all}\ Z > 0, \ \lim_{Z \to \infty} \frac{P(Z)}{Z^{5/3}} = P_\infty > 0.
\eF

\item In agreement with Gibbs' relation (\ref{i6}), we take
\bFormula{H4}
e_M (\vr, \vt) = \frac{3}{2} \vt \left( \frac{ \vt^{3/2} }{\vr} \right) P \left( \frac{\vr}{\vt^{3/2}} \right),
\eF
and
\bFormula{H5}
s_M (\vr, \vt) = S \left( \frac{\vr}{\vt^{3/2}} \right),
\eF
where
\bFormula{H6}
S'(Z) = - \frac{3}{2} \frac{ \frac{5}{3} P(Z) - P'(Z) Z }{Z^2} < 0.
\eF
In addition, we require
\bFormula{H7}
\lim_{Z \to \infty} S(Z) = 0.
\eF
\item The viscosity coefficients in (\ref{e5}) are continuously differentiable functions of $\vt \in [0, \infty)$ satisfying
\bFormula{H8}
| \mu'(\vt) | \leq c, \ \underline{\mu} (1 + \vt) \leq \mu(\vt) ,\ 0 \leq \eta(\vt) \leq \Ov{\eta}(1 + \vt)
\  \mbox{for all}\ \vt \geq 0
\eF
for certain constants $\underline{\mu} > 0$, $\Ov{\eta} > 0$.
\item The heat conductivity coefficient in Fourier's law (\ref{e6}) satisfies
\bFormula{H9}
\kappa \in C^1[0, \infty),\ \underline{\kappa} (1 + \vt^3) \leq \kappa(\vt) \leq \Ov{\kappa}(1 + \vt^3) \ \mbox{for all}\ \vt \geq 0
\eF
for certain constants $\underline{\kappa} > 0$, $\Ov{\kappa} > 0$.
\end{itemize}

The interested reader may consult \cite[Chapter 2]{FEINOV} for the physical background as well as possible relaxation of these restrictions. The important fact is
that imposing the hypotheses (\ref{H1}--\ref{H7}), together with (\ref{e7}--\ref{e8a}), render the Navier-Stokes-Fourier system mathematically tractable, at least in the context of weak solutions. Specifically, as shown in \cite[Chapter 3, Theorem 3.1]{FEINOV}, the problem (\ref{e1}--\ref{e21}) admits a global-in-time weak solution for any choice of finite energy initial data $[\vr_0, \vt_0, \vu_0]$ satisfying the natural constraint $\vr_0 > 0$, $\vt_0 > 0$.

\bRemark{H1}

It is worth noting that the problem considered in \cite[Chapter 3]{FEINOV} does not involve the damping term $-\lambda \vu$ in the momentum equation. However, the existence
proof  can be easily adapted to the present setting yielding the total energy inequality (\ref{e4}) in place of the equality asserted by \cite[Theorem 3.1]{FEINOV}.

\eR

\section{Main result}
\label{m}

Following \cite{FeiNov10}, we introduce the ballistic free energy
\[
H_\Theta(\vr, \vt) = \vr \Big( e(\vr, \vt) - \Theta s(\vr, \vt) \Big),
\]
and the relative energy functional
\bFormula{rr1}
\mathcal{E} \left( \vr, \vt, \vu \Big| r, \Theta, \vc{U} \right) =
\intO{ \left[ \frac{1}{2} \vr |\vu - \vc{U} |^2 + H_\Theta (\vr, \vt) - \frac{\partial H_\Theta(r, \Theta) }{\partial \vr} (\vr - r) -
H_\Theta (r, \Theta) \right] }.
\eF

As shown in \cite[Chapter 5, Lemma 5.1]{FEINOV}, the functional $\mathcal{E}$ can be viewed as a kind of ``distance'' between the quantities
$[\vr, \vt, \vu]$ and $[r, \Theta, \vc{U}]$. Specifically, for any compact set $K \subset (0, \infty)^2$, there exists a positive constant
$c(K)$, depending solely on the structural properties of the thermodynamic functions stated in (\ref{H1}--\ref{H7}) such that
\bFormula{RrR1}
\mathcal{E} \left( \vr, \vt, \vu \Big| r, \Theta, \vc{U} \right) \geq c(K) \left\{ \begin{array}{l}  |\vr - r|^2 + |\vt - \Theta|^2 + |\vu - \vc{U}|^2
\ \mbox{if}\ [\vr, \vt] \in K, \ [r, \Theta] \in K \\ \\
1 + \vr |\vu - \vc{U}|^2 + \vr e(\vr, \vt) + \vr |s(\vr, \vt)|  \ \mbox{if} \ [\vr, \vt] \in (0, \infty)^2 \setminus K, \ [r, \Theta] \in K.
\end{array} \right.
\eF

In order to state our main result, we need strong solutions to the Euler system (\ref{i1}--\ref{i3}) supplemented with the boundary condition
\bFormula{m1}
\vu \cdot \vc{n}|_{\partial \Omega} = 0.
\eF
Note that the total energy balance (\ref{i3}) can be equivalently reformulated as the entropy balance equation
\bFormula{m2}
\partial_t (\vr s_M(\vr, \vt)) + \Div (\vr s_M(\vr, \vt) \vu) = 0,
\eF
or the thermal energy balance
\bFormula{m3}
c_v (\vr, \vt) \left( \partial_t (\vr \vt) + \Div (\vr \vt \vu) \right) + \vt \frac{\partial p_M (\vr, \vt)}{\partial \vt} \Div \vu = 0,
\ \mbox{with}\
c_v (\vr, \vt) = \frac{\partial e_M(\vr, \vt)}{\partial \vt},
\eF
as long as the solution of the Euler system remains smooth.

A suitable existence result for the Euler system with the slip boundary condition (\ref{m1}) was obtained by Schochet \cite[Theorem 1]{SCHO1}. It asserts the local-in-time existence of a \emph{classical} solution $[\vr_E, \vt_E, \vu_E]$ of the Euler system (\ref{i1}), (\ref{i2}), (\ref{m1}), (\ref{m2}) if:
\begin{itemize}
\item $\Omega \subset R^3$ is a bounded domain with a sufficiently smooth boundary, say $\partial \Omega$ of class $C^\infty$;
\item the initial data $[\vr_{0,E}, \vt_{0,E}, \vu_{0,E}]$ satisfy
\bFormula{m4}
\vr_{0, E}, \vt_{0, E} \in W^{3,2}(\Omega), \ \vu_{0, E} \in W^{3,2}(\Omega; R^3),\
\vr_{0, E}, \ \vt_{0,E} > 0  \ \mbox{in}\ \Ov{\Omega};
\eF
\item the compatibility conditions
\bFormula{m5}
\partial^k_t \vu_{0,E} \cdot \vc{n}|_{\partial \Omega} = 0
\eF
hold for $k=0,1,2$.

\end{itemize}

We are ready to state the main result of this paper.

\bTheorem{m1}
Let $\Omega \subset R^3$ be a bounded domain with smooth boundary. Suppose that the thermodynamic functions $p$, $e$, and $s$ are given by (\ref{e7}--\ref{e8a}), where $p_M$, $e_M$, and $s_M$ satisfy (\ref{H1}--\ref{H7}). Let the transport coefficients $\mu$, $\eta$ and $\lambda$ obey (\ref{H8}), (\ref{H9}). Let
$[\vr_E, \vt_E, \vu_E]$ be the classical solution of the Euler system (\ref{i1}--\ref{i3}), (\ref{m1}) in a time interval $(0,T)$, originating from the initial data $[\vr_{0,E}, \vt_{0,E}, \vu_{0,E}]$ satisfying (\ref{m4}), (\ref{m5}). Finally, let $[\vr, \vt, \vu]$ be a weak solution of the Navier-Stokes-Fourier system (\ref{e1}--\ref{e6}), (\ref{e20}), (\ref{e21}), where the initial data $[\vr_0, \vt_0, \vu_0]$ satisfy
\bFormula{m6}
\vr_0, \vt_0 > 0 \ \mbox{a.a. in} \ \Omega,\ \intO{ \vr_0 } \geq M, \ \|\vr_0 \|_{L^\infty(\Omega)} + \|\vt_0 \|_{L^\infty(\Omega)} + \|\vu_0 \|_{L^\infty(\Omega;R^3)} \leq D,
\eF
and where the scaling parameters $a$, $\nu$, $\omega$, and $\lambda$ are positive numbers.

Then
\bFormula{konec}
\mathcal{E} \left( \vr, \vt, \vu \Big| \vrE, \vtE, \vuE \right)(\tau)
\eF
\[
\aleq c(T,M,D) \left( \mathcal{E} \left( \vr_0, \vt_0, \vu_0 \Big| \vr_{0,E}, \vt_{0,E}, \vu_{0,E} \right)
+ \max \left\{a, \nu, \omega, \lambda,
\frac{\nu}{\sqrt{a}}, \frac{\omega}{a}, \left( \frac{a}{ \sqrt{\nu^3
\lambda} } \right)^{1/3} \right\} \right)
\]
for a.a. $\tau \in (0,T)$.
\eT

In view of the coercivity properties of the relative energy specified in (\ref{RrR1}), Theorem \ref{Tm1} provides an immediate corollary.

\bCorollary{m1}
Under the hypotheses of Theorem \ref{Tm1} suppose that
\bFormula{m7}
a, \nu, \omega, \lambda \to 0,\ \mbox{and}\
\frac{\omega}{a} \to 0, \ \frac{\nu}{\sqrt{a}} \to 0,\ \frac{a}{ \sqrt{ \nu^3 \lambda } } \to 0.
\eF

Then
\bFormula{m8}
{\rm ess} \sup_{\tau \in (0,T)} \intO{ \left[
\vr |\vu - \vu_E |^2 + | \vr - \vr_E |^{5/3} + \vr |\vt - \vt_E | \right] }
\eF
\[
\leq c(T,D,M) \Lambda \left(a, \nu, \omega, \lambda, \| \vr_0 - \vr_{0,E} \|_{L^\infty(\Omega)}, \| \vt_0 - \vt_{0,E} \|_{L^\infty(\Omega)},
\| \vu_0 - \vu_{0,E} \|_{L^\infty(\Omega;R^3)} \right),
\]
where $\Lambda$ is an explicitly computable function of its arguments,
\[
\Lambda \left(a, \nu, \omega, \lambda, \| \vr_0 - \vr_{0,E} \|_{L^\infty(\Omega)}, \| \vt_0 - \vt_{0,E} \|_{L^\infty(\Omega)},
\| \vu_0 - \vu_{0,E} \|_{L^\infty(\Omega;R^3)} \right) \to 0
\]
provided $a, \nu, \omega, \lambda$ satisfy (\ref{m7}), and
\[
\| \vr_0 - \vr_{0,E} \|_{L^\infty(\Omega)}, \| \vt_0 - \vt_{0,E} \|_{L^\infty(\Omega)},
\| \vu_0 - \vu_{0,E} \|_{L^\infty(\Omega;R^3)} \to 0.
\]

\eC

\bRemark{m1}

As already mentioned above, the convergence stated in Corollary \ref{Cm1} is path dependent, the parameters $a, \nu, \omega, \lambda$ are interrelated through
(\ref{m7}). It is easy to see that (\ref{m7}) holds provided, for instance,
\[
a \to 0,\ \nu = a^\alpha, \ \omega = a^\beta, \ \lambda = a^\gamma, \ \mbox{where}\
\beta > 1, \  \frac{1}{2} < \alpha < \frac{2}{3}, \ 0 < \gamma < 1  - \frac{3}{2} \alpha.
\]

\eR

\bRemark{m2}

The specific form of the function $\Lambda$ is given in terms of the relative energy functional $\mathcal{E}$ introduced at the beginning of this section.
The constants in (\ref{konec}) and (\ref{m8}) depend also on the properties of the limit solution $[\vrE, \vtE, \vuE]$.

\eR

The rest of the paper is devoted to the proof of Theorem \ref{Tm1}. Possible extensions are discussed in Section \ref{F}.

\section{Relative energy}
\label{r}

As shown in \cite{FeiNov10}, \emph{any} weak solution $[\vr, \vt, \vu]$ to the Navier-Stokes-Fourier system with the boundary conditions (\ref{e20}), (\ref{e21}) satisfies the relative energy inequality
\bFormula{r1}
\left[ \mathcal{E} \left( \vr , \vt , \vu \Big| r, \Theta, \vc{U} \right) \right]_{t = 0}^{t = \tau}
\eF
\[
 +
\int_0^\tau \intO{ \frac{\Theta}{\vt} \left(
\tn{S}(\vt, \Grad \vu): \Grad \vu - \frac{ \vc{q}(\vt, \Grad \vt) \cdot \Grad \vt }{\vt} \right) } \ \dt + \lambda \int_0^\tau \intO{ |\vu|^2 } \ \dt
\]
\[
\leq
\int_0^\tau \intO{ \vr (\vu - \vc{U}) \cdot \Grad \vc{U} \cdot (\vc{U} - \vu)} \ \dt
\]
\[
+ \int_0^\tau \intO{  \tn{S}(\vt, \Grad \vu) : \Grad \vc{U} } \ \dt - \int_0^\tau \intO{
\frac{ \vc{q}(\vt, \Grad \vt) }{\vt} \cdot \Grad \Theta  } \ \dt +  \lambda \int_0^\tau \intO{ \vu \cdot \vc{U} } \ \dt
\]
\[
+ \int_0^\tau \intO{ \vr \Big( s(\vr, \vt) - s(r, \Theta) \Big) \Big( \vc{U} - \vu \Big)
\cdot \Grad \Theta } \ \dt
\]
\[
+ \int_0^\tau \intO{ \vr \Big(  \partial_t \vc{U} +  \vc{U} \cdot \Grad \vc{U} \Big) \cdot (\vc{U} - \vu) } \ \dt
- \int_0^\tau \intO{ p(\vr, \vt) \Div \vc{U}  } \ \dt
\]
\[
- \int_0^\tau \intO{ \left( \vr \Big( s (\vr, \vt) - s(r, \Theta) \Big) \partial_t \Theta
+ \vr \Big( s(\vr, \vt) - s(r, \Theta) \Big) \vc{U} \cdot \Grad \Theta \right) } \ \dt
\]
\[
+ \int_0^\tau \intO{ \left( \left( 1 - \frac{\vr}{r} \right) \partial_t p(r, \Theta) -
\frac{\vr}{r} \vu \cdot \Grad p(r, \Theta) \right) } \ \dt
\]
for any trio of (smooth) test functions $[r, \Theta, \vc{U}]$ such that
\bFormula{rr2}
r, \ \Theta > 0  \ \mbox{in}\ \Ov{\Omega}, \ \vc{U} \cdot \vc{n}|_{\partial \Omega} = 0.
\eF

\bRemark{rr1}

Following Lions' idea \cite{LI} proposed for the incompressible Euler system, we may define \emph{dissipative solutions} to the Navier-Stokes-Fourier
system requiring solely certain regularity of $[\vr, \vt, \vu]$ and the relative energy inequality (\ref{r1}) to be satisfied for any trio
$[r, \Theta, \vc{U}]$ as in (\ref{rr2}). Such an approach was used by Jessl\' e, Jin, Novotn\' y \cite{JeJiNo} to attack problems on unbounded spatial domains.

\eR

\section{Energy estimates}
\label{E}

In what follows, we use the notation
\[
a \aleq b \ \mbox{if} \ a \leq cb
\]
where $c > 0$ is a constant independent of the scaling parameters $a, \nu, \omega, \lambda$, and of $D$ and $T$. Similarly, we define $a \ageq b$ and
$a \approx b$.

We start with some auxiliary estimates that follow directly from the structural hypotheses imposed on the functions $p$, $e$, and $s$, see
\cite[Chapter 3]{FEINOV} for the proofs:
\bFormula{e9}
0 \leq \vr s_M(\vr, \vt) \aleq \vr \left( 1 + |\log(\vr)| + [\log(\vt)]^+ \right),
\eF
\bFormula{E2}
\vr e_M (\vr, \vt) \ageq \vr \vt + \vr^{5/3}.
\eF

Taking $r>0$, $\Theta>0$ constant and $\vc{U} = 0$ in the relative energy inequality (\ref{r1}), and keeping (\ref{e9}), (\ref{E2}) in mind, we deduce the standard energy estimates:
\bFormula{V1a}
\left[
\begin{array}{c}
{\rm ess} \sup_{t \in (0,T)} \intO{ \vr |\vu|^2(t, \cdot) } \leq c(D),\\ \\
{\rm ess} \sup_{t \in (0,T)} \intO{ \vr^{5/3}(t, \cdot) } \leq c(D), \\ \\
{\rm ess} \sup_{t \in (0,T)} \intO{ \vr \vt (t, \cdot)} \leq c(D),\\ \\
a \ {\rm ess} \sup_{t \in (0,T)} \intO{ \vt^4 (t, \cdot) } \leq c(D);
\end{array}
\right]
\eF
together with the estimates following from dissipation and hypotheses (\ref{H8}), (\ref{H9}):
\bFormula{V1}
\left[
\begin{array}{c}
\nu \int_0^T \intO{ \left|\Grad \vu + \Grad^t \vu - \frac{2}{3} \Div \vu \tn{I} \right|^2 } \ \dt \leq c(D),\
\lambda \int_0^T \intO{ |\vu|^2 } \ \dt \leq c(D),\\ \\
\omega \int_0^T \intO{ \left[ |\Grad \vt |^2 + |\log(\vt)|^2 \right] } \ \dt \leq c(D).
\end{array}
\right]
\eF

Now, by virtue of hypothesis (\ref{m6}) and the bounds (\ref{V1a} $)_{1,2}$, (\ref{V1}$)_1$, we may us a generalized version of Korn's inequality
\cite[Theorem 10.17]{FEINOV} to obtain
\bFormula{V1b}
\| \sqrt{\nu} \vu \|_{L^2(0,T; W^{1,2}(\Omega; R^3)} \leq c(M,D);
\eF
whence,
my means of the standard embedding relations for Sobolev functions,
\bFormula{V2}
\| \sqrt{\nu} \vu \|_{L^2(0,T; L^6(\Omega; R^3))} \leq c(M,D).
\eF

Moreover, in view of the interpolation inequality
\[
\| \vu \|_{L^4(\Omega; R^3)} \leq \| \vu \|_{L^6(\Omega; R^3)}^{3/4} \| \vu \|_{L^2(\Omega; R^3))}^{1/4},
\]
we obtain
\bFormula{V2a}
\left\| \nu^{3/8} \lambda^{1/8} \vu \right\|_{L^4(\Omega; R^3)} \leq \| \sqrt{\nu} \vu \|_{L^6(\Omega; R^3)}^{3/4} \| \sqrt{\lambda} \vu \|_{L^2(\Omega; R^3)}^{1/4} \aleq \left( \| \sqrt{\nu} \vu \|_{L^6(\Omega; R^3)} + \| \sqrt{\lambda} \vu \|_{L^2(\Omega; R^3)} \right),
\eF
therefore, by virtue of (\ref{V1} $)_1$, (\ref{V2}),
\bFormula{V3a}
\left\| \nu^{3/8} \lambda^{1/8} \vu \right\|_{L^2(0,T; L^4(\Omega; R^3))} \leq c(M,D).
\eF

\section{Convergence}
\label{c}

Our goal is to show Theorem \ref{Tm1}. The obvious idea is to take
\bFormula{cc1-}
r = \vr_E, \ \Theta = \vt_E, \ \vc{U} = \vu_E
\eF
as test functions in the relative energy inequality (\ref{r1}). We emphasize again that such a step is conditioned by our choice of the slip boundary conditions for the velocity field in both the primitive and the target system.

Next, we fix positive constants $\underline{\vr}$, $\Ov{\vr}$, $\underline{\vt}$, $\Ov{\vt}$ in such a way that
\[
0 < \underline{\vr} < \inf_{t \in [0,T], x \in \Ov{\Omega}} \vr_E(t,x) \leq \sup_{t \in [0,T], x \in \Ov{\Omega}} \vr_E(t,x) < \Ov{\vr},
\]
\[
0 < \underline{\vt} < \inf_{t \in [0,T], x \in \Ov{\Omega}} \vt_E(t,x) \leq \sup_{t \in [0,T], x \in \Ov{\Omega}} \vt_E(t,x) < \Ov{\vt}.
\]
Following \cite{FEINOV} we introduce a decomposition of a measurable function $F$ into its essential and residual part, specifically,
\[
F = \ess{F} + \res{F}, \ \ess{F} = \Phi (\vr_E, \vt_E) F, \ \res{F} = (1 - \Phi(\vr_E, \vt_E)) F,
\]
where
\[
\Phi \in \DC((0, \infty)^2), \ 0 \leq \Phi \leq 1, \ \Phi \equiv 1 \ \mbox{on the rectangle}\ [\underline{\vr}, \Ov{\vr}] \times [\underline \vt, \Ov{\vt}].
\]

Now, it follows from hypotheses (\ref{H2}), (\ref{H6}), and the estimates stated in (\ref{RrR1}) that
\bFormula{cc1}
\left\| \ess{\vr - \vrE} \right\|_{L^2(\Omega)}^2 + \left\| \ess{\vt - \vtE} \right\|_{L^2(\Omega)}^2 +
\left\| \ess{\vu - \vuE} \right\|_{L^2(\Omega; R^3)}^2 \aleq \mathcal{E} \left( \vr, \vt, \vu \Big| \vrE, \vtE, \vuE \right) ,
\eF
and, by virtue of (\ref{e9}), (\ref{E2}),
\bFormula{cc2}
\intO{ \vr |\vu - \vuE|^2 } +
\intO{ \res{ 1 + \vr^{5/3} + \vr \vt + a \vt^4 } } \aleq  \mathcal{E} \left( \vr, \vt, \vu \Big| \vrE, \vtE, \vuE \right) .
\eF

Using the ansatz (\ref{cc1-}) we examine separately all integrals on the right-hand side of (\ref{r1}). As, obviously,
\[
\int_0^\tau \intO{ \vr (\vu - \vuE) \cdot \Grad \vuE \cdot (\vuE - \vu)} \ \dt \leq \int_0^\tau \mathcal{E} \left( \vr, \vt, \vu \Big| \vrE, \vtE, \vuE \right),
\]
the relative energy inequality (\ref{r1}) can be written as
\bFormula{vhr1}
\left[ \mathcal{E} \left( \vr , \vt , \vu \Big| \vrE, \vtE, \vuE \right) \right]_{t = 0}^{t = \tau}
\eF
\[
 +
\int_0^\tau \intO{ \frac{\vtE}{\vt} \left(
\tn{S}(\vt, \Grad \vu): \Grad \vu - \frac{ \vc{q}(\vt, \Grad \vt) \cdot \Grad \vt }{\vt} \right) } \ \dt + \lambda \int_0^\tau \intO{ |\vu|^2 } \ \dt
\]
\[
\aleq
\int_0^\tau \mathcal{E} \left( \vr , \vt , \vu \Big| \vrE, \vtE, \vuE \right) \dt
\]
\[
+ \int_0^\tau \intO{  \tn{S}(\vt, \Grad \vu) : \Grad \vuE } \ \dt - \int_0^\tau \intO{
\frac{ \vc{q}(\vt, \Grad \vt) }{\vt} \cdot \Grad \vtE  } \ \dt +  \lambda \int_0^\tau \intO{ \vu \cdot \vuE } \ \dt
\]
\[
+ \int_0^\tau \intO{ \vr \Big( s(\vr, \vt) - s(\vrE, \vtE) \Big) \Big( \vuE - \vu \Big)
\cdot \Grad \vtE } \ \dt
\]
\[
+ \int_0^\tau \intO{ \vr \Big(  \partial_t \vuE +  \vuE \cdot \Grad \vuE \Big) \cdot (\vuE - \vu) } \ \dt
- \int_0^\tau \intO{ p(\vr, \vt) \Div \vuE  } \ \dt
\]
\[
- \int_0^\tau \intO{ \left( \vr \Big( s (\vr, \vt) - s(\vrE, \vtE) \Big) \Big( \partial_t \vtE
+ \vuE \cdot \Grad \vtE \right) } \ \dt
\]
\[
+ \int_0^\tau \intO{ \left( \left( 1 - \frac{\vr}{\vrE} \right) \partial_t p(\vrE, \vtE) -
\frac{\vr}{\vrE} \vu \cdot \Grad p(\vrE, \vtE) \right) } \ \dt.
\]

\subsection{Integrals depending on viscosity and the heat flux}

On one hand, we have
\bFormula{vh1}
\frac{\vtE}{\vt} \tn{S} (\vt, \Grad \vu) : \Grad \vu \ageq \nu \left( \frac{\mu(\vt)}{\vt} \left| \Grad \vu + \Grad^t \vu - \frac{2}{3} \Div \vu \tn{I} \right|^2 +
\frac{\eta(\vt)}{\vt} | \Div \vu |^2 \right),
\eF
while, on the other hand,
\bFormula{vh2}
\tn{S}(\vt, \Grad \vu) : \Grad \vuE
\eF
\[
= \nu \left( \frac{\mu(\vt)}{\vt} \right)^{1/2} \left( \Grad \vu + \Grad^t \vu - \frac{2}{3} \Div \vu \tn{I} \right) : \left( {\mu(\vt)}{\vt} \right)^{1/2} \Grad \vuE
+ \nu \left( \frac{\eta(\vt)}{\vt} \right)^{1/2} \Div \vu  \left( {\eta(\vt)}{\vt} \right)^{1/2} \Div \vuE
\]
\[
\leq \delta \nu  \left( \frac{\mu(\vt)}{\vt} \left| \Grad \vu + \Grad^t \vu - \frac{2}{3} \Div \vu \tn{I} \right|^2 +
\frac{\eta(\vt)}{\vt} | \Div \vu |^2 \right)
+ \nu c(\delta) \vt \Big( \mu(\vt) + \eta(\vt) \Big) |\Grad \vuE|^2 .
\]
\[
\leq \frac{1}{2} \frac{\vtE}{\vt} \tn{S} (\vt, \Grad \vu) : \Grad \vu + \nu c(\delta) \vt \Big( \mu(\vt) + \eta(\vt) \Big) |\Grad \vuE |^2
\]
provided the parameter $\delta > 0$ has been taken small enough. Moreover, in view of hypotheses (\ref{H8}), (\ref{m7}) and the energy bound
(\ref{V1a} $)_4$,
\bFormula{vh2a}
\nu \intO{ \vt \Big( \mu(\vt) + \eta(\vt) \Big) |\Grad \vuE |^2 }
\eF
\[
\aleq \mathcal{O}(\nu) + \intO{ \nu \vt^2 }
\aleq \mathcal{O}(\nu) + \left( \intO{ \nu^2 \vt^4 } \right)^{1/2} \aleq \mathcal{O}(\nu) + \mathcal{O} \left( \frac{\nu}{\sqrt{a}} \right)
\to 0.
\]

Similarly,
\bFormula{vh3}
- \frac{\vtE}{\vt} \frac{\vc{q} (\vt, \Grad \vt) \cdot \Grad \vt}{\vt} \ageq \omega \frac{\kappa(\vt)}{\vt^2} |\Grad \vt|^2 \ageq
\omega \left( |\Grad \vt |^2 + |\Grad \log (\vt) |^2 \right),
\eF
while
\bFormula{vh4}
- \frac{ \vc{q}(\vt, \Grad \vt) }{\vt} \cdot \Grad \vtE = \omega \left( \frac{\kappa (\vt) }{\vt^2} \right)^{1/2} \Grad \vt \cdot
\left( \kappa(\vt) \right)^{1/2} \Grad \vtE
\eF
\[
\aleq \delta \omega \frac{\kappa(\vt)}{\vt^2} |\Grad \vt|^2 + \omega c(\delta) \kappa(\vt) |\Grad \vtE|^2
\leq - \frac{1}{2} \frac{\vtE}{\vt} \frac{\vc{q} (\vt, \Grad \vt) \cdot \Grad \vt}{\vt} + \omega c(\delta) \kappa(\vt) |\Grad \vtE|^2,
\]
where, furthermore,
\bFormula{vh4a}
\omega \intO{ \kappa(\vt) |\Grad \vtE|^2 } \aleq \mathcal{O}(\omega) + \frac{\omega}{a} \intO{ a \vt^4 } \aleq
\mathcal{O}(\omega) + \mathcal{O} \left( \frac{\omega}{a} \right).
\eF

Finally, seeing that
\bFormula{vh4b}
\lambda \intO{ \vu \cdot \vuE } \leq \frac{\lambda}{2} \intO{ |\vu|^2 } + \mathcal{O} (\lambda),
\eF
we may combine (\ref{vh1}--\ref{vh4b}) to reduce (\ref{vhr1}) to
\bFormula{vhr2}
\left[ \mathcal{E} \left( \vr , \vt , \vu \Big| \vrE, \vtE, \vuE \right) \right]_{t = 0}^{t = \tau}
\aleq
\int_0^\tau \mathcal{E} \left( \vr , \vt , \vu \Big| \vrE, \vtE, \vuE \right) \dt + c(E,M) \mathcal{O} \left( \max \left\{ \nu, \omega, \lambda,
\frac{\nu}{\sqrt{a}}, \frac{\omega}{a} \right\} \right)
\eF
\[
+ \int_0^\tau \intO{ \vr \Big( s(\vr, \vt) - s(\vrE, \vtE) \Big) \Big( \vuE - \vu \Big)
\cdot \Grad \vtE } \ \dt
\]
\[
+ \int_0^\tau \intO{ \vr \Big(  \partial_t \vuE +  \vuE \cdot \Grad \vuE \Big) \cdot (\vuE - \vu) } \ \dt
- \int_0^\tau \intO{ p(\vr, \vt) \Div \vuE  } \ \dt
\]
\[
- \int_0^\tau \intO{ \left( \vr \Big( s (\vr, \vt) - s(\vrE, \vtE) \Big) \Big( \partial_t \vtE
+ \vuE \cdot \Grad \vtE \right) } \ \dt
\]
\[
+ \int_0^\tau \intO{ \left( \left( 1 - \frac{\vr}{\vrE} \right) \partial_t p(\vrE, \vtE) -
\frac{\vr}{\vrE} \vu \cdot \Grad p(\vrE, \vtE) \right) } \ \dt.
\]

\subsection{Entropy integral}

The first integral on the right hand side of (\ref{vhr2}) is the most difficult one to handle. It can be written in the form
\bFormula{vh5}
\int_0^\tau \intO{ \vr \Big( s(\vr, \vt) - s(\vrE, \vtE) \Big) \Big( \vuE - \vu \Big)
\cdot \Grad \vtE } \ \dt
\eF
\[
= \int_0^\tau \intO{ \vr \ess{ s(\vr, \vt) - s(\vrE, \vtE) } \Big( \vuE - \vu \Big)
\cdot \Grad \vtE } \ \dt
\]
\[
+ \int_0^\tau \intO{ \vr \res{ s(\vr, \vt) - s(\vrE, \vtE) } \Big( \vuE - \vu \Big)
\cdot \Grad \vtE } \ \dt,
\]
where, in view of (\ref{cc1}),
\[
\left|
\int_0^\tau \intO{ \vr \ess{ s(\vr, \vt) - s(\vrE, \vtE) } \Big( \vuE - \vu \Big)
\cdot \Grad \vtE } \ \dt \right| \aleq \int_0^\tau \mathcal{E} \left( \vr, \vt, \vu \Big| \ \vrE, \vtE, \vuE \right) \ \dt.
\]

Next, in accordance with (\ref{e9}),
\bFormula{vh6}
\left| \int_0^\tau \intO{ \vr \res{ s(\vr, \vt) - s(\vrE, \vtE) } \Big( \vuE - \vu \Big)
\cdot \vtE } \ \dt \right|
\eF
\[
\aleq \int_0^\tau \intO{ a \res{ \vt^3 } \left( |\vu| + 1 \right) } \ \dt + \int_0^\tau \intO{ \sqrt{\vr} \res{ 1 + |\log(\vr)| + [\log(\vt)]^+
} \sqrt{\vr} |\vu - \vuE| } \ \dt,
\]
where, by means of (\ref{cc2}),
\bFormula{vh6a}
\int_0^\tau \intO{ \sqrt{\vr} \res{ 1 + |\log(\vr)| + [\log(\vt)]^+
} \sqrt{\vr} |\vu - \vuE| } \ \dt
\eF
\[
\aleq \int_0^\tau \intO{ \res{ \vr + \vr^{5/3} + \vr \vt } } \ \dt +
\int_0^\tau \intO{ \vr |\vu - \vuE|^2 } \ \dt
\aleq \int_0^\tau \mathcal{E} \left( \vr, \vt, \vu \Big|  \vrE, \vtE, \vuE \right) \ \dt.
\]

Finally, by H\" older's inequality and (\ref{cc2}),
\[
\intO{ a \res{ \vt^3 }  ( |\vu| + 1) } \leq \left( \intO{ \frac{a^{4/3}}{ \nu^{1/2} \lambda^{1/6} } \vt^4 } \right)^{3/4} \left\| \nu^{3/8} \lambda^{1/8} \vu \right\|_{L^4(\Omega; R^3)} + \mathcal{E} \left( \vr, \vt, \vu \Big|  \vrE, \vtE, \vuE \right),
\]
where, furthermore,
\[
\intO{ \frac{a^{4/3}}{ \nu^{1/2} \lambda^{1/6} } \vt^4 }  = \frac{a^{1/3}}{\nu^{1/2} \lambda^{1/6}} \intO{ a \vt^4 } \aleq c(E) \left( \frac{a}{ \sqrt{\nu^3
\lambda} } \right)^{1/3} \to 0.
\]

Summing up the previous estimates, we may rewrite (\ref{vhr2}) as
\bFormula{vhr3}
\left[ \mathcal{E} \left( \vr , \vt , \vu \Big| \vrE, \vtE, \vuE \right) \right]_{t = 0}^{t = \tau}
\eF
\[
\aleq
\int_0^\tau \mathcal{E} \left( \vr , \vt , \vu \Big| \vrE, \vtE, \vuE \right) \dt + c(E,M) \mathcal{O} \left( \max \left\{ \nu, \omega, \lambda,
\frac{\nu}{\sqrt{a}}, \frac{\omega}{a}, \left( \frac{a}{ \sqrt{\nu^3
\lambda} } \right)^{1/3} \right\} \right)
\]
\[
+ \int_0^\tau \intO{ \vr \Big(  \partial_t \vuE +  \vuE \cdot \Grad \vuE \Big) \cdot (\vuE - \vu) } \ \dt
- \int_0^\tau \intO{ p(\vr, \vt) \Div \vuE  } \ \dt
\]
\[
- \int_0^\tau \intO{ \left( \vr \Big( s (\vr, \vt) - s(\vrE, \vtE) \Big) \Big( \partial_t \vtE
+ \vuE \cdot \Grad \vtE \right) } \ \dt
\]
\[
+ \int_0^\tau \intO{ \left( \left( 1 - \frac{\vr}{\vrE} \right) \partial_t p(\vrE, \vtE) -
\frac{\vr}{\vrE} \vu \cdot \Grad p(\vrE, \vtE) \right) } \ \dt.
\]

\subsection{Remaining integrals}

Keeping in mind that $[\vrE, \vtE, \vuE]$ is a smooth solution of the Euler system, we may handle
the remaining integrals on the right-hand side of (\ref{vhr3}) as follows:
\bFormula{oi1}
\int_0^\tau \intO{ \vr \Big(  \partial_t \vuE +  \vuE \cdot \Grad \vuE \Big) \cdot (\vuE - \vu) } \ \dt
- \int_0^\tau \intO{ p(\vr, \vt) \Div \vuE   } \ \dt
\eF
\[
- \int_0^\tau \intO{  \vr \Big( s (\vr, \vt) - s(\vrE, \vtE) \Big) \Big[ \partial_t \vtE
+ \vuE \cdot \Grad \vtE \Big] } \ \dt
\]
\[
+ \int_0^\tau \intO{ \left( \left( 1 - \frac{\vr}{\vrE} \right) \partial_t p(\vrE, \vtE) -
\frac{\vr}{\vrE} \vu \cdot \Grad p(\vrE, \vtE) \right) } \ \dt
\]
\[
= - \int_0^\tau \intO{ \frac{\vr}{\vrE} \Grad p_M (\vrE, \vtE) \cdot (\vuE - \vu) } \ \dt
- \int_0^\tau \intO{ p(\vr, \vt) \Div \vuE   } \ \dt
\]
\[
- \int_0^\tau \intO{  \vr \Big( s (\vr, \vt) - s(\vrE, \vtE) \Big) \Big[ \partial_t \vtE
+ \vuE \cdot \Grad \vtE \Big] } \ \dt
\]
\[
+ \int_0^\tau \intO{ \left( \left( 1 - \frac{\vr}{\vrE} \right) \partial_t p_M(\vrE, \vtE) -
\frac{\vr}{\vrE} \vu \cdot \Grad p_M(\vrE, \vtE) \right) } \ \dt + \mathcal{O}(a)
\]
\[
=
\int_0^\tau \intO{ \Big[ p_M(\vrE, \vtE)  - p(\vr, \vt) \Big] \Div \vuE   } \ \dt
\]
\[
- \int_0^\tau \intO{  \vr \Big( s (\vr, \vt) - s(\vre, \vtE) \Big) \Big[ \partial_t \vtE
+ \vuE \cdot \Grad \vtE \Big] } \ \dt
\]
\[
+ \int_0^\tau \intO{  \left( 1 - \frac{\vr}{\vrE} \right) \left[ \partial_t p_M(\vrE, \vtE) +
\vuE \cdot \Grad p_M(\vrE, \vtE) \right] } \ \dt + \mathcal{O}(a)
\]
\[
\aleq \int_0^\tau \intO{ \ess{ p_M(\vrE, \vtE)  - p_M(\vr, \vt) } \Div \vuE   } \ \dt
\]
\[
- \int_0^\tau \intO{  \vr \Big( s_M (\vr, \vt) - s_M (\vrE, \vtE) \Big) \ess{ \partial_t \vtE
+ \vuE \cdot \Grad \vtE } } \ \dt
\]
\[
+ \int_0^\tau \intO{  \left( 1 - \frac{\vr}{\vrE} \right) \ess{ \partial_t p_M(\vrE, \vtE) +
\vuE \cdot \Grad p_M(\vrE, \vtE) } } \ \dt
\]
\[
+ \int_0^\tau \intO{ \mathcal{E} \left( \vr, \vt, \vu \Big| \vrE, \vtE, \vuE \right) } + \mathcal{O}(a).
\]

Now, using the equations satisfied by $[\vrE, \vtE, \vuE]$, together with Gibbs' relation (\ref{i6}), we may write
\[
\partial_t p_M(\vrE, \vtE) + \vuE \cdot \Grad p_M(\vrE, \vtE)
\]
\[
= \frac{\partial p_M (\vrE, \vtE)}{\partial \vr} \left( \partial_t \vrE + \vuE \cdot \Grad \vrE \right)
+\frac{\partial p_M (\vrE, \vtE)}{\partial \vt} \left( \partial_t \vtE + \vuE \cdot \Grad \vtE \right)
\]
\[
=- \vrE \frac{\partial p_M (\vrE, \vtE)}{\partial \vr} \Div \vuE
- \vrE^2 \frac{\partial s_M (\vrE, \vtE)}{\partial \vr} \left( \partial_t \vtE + \vuE \cdot \Grad \vtE \right);
\]
whence (\ref{oi1}) gives rise to
\bFormula{oi2}
\int_0^\tau \intO{ \vr \Big(  \partial_t \vuE +  \vuE \cdot \Grad \vuE \Big) \cdot (\vuE - \vu) } \ \dt
- \int_0^\tau \intO{ p(\vr, \vt) \Div \vuE   } \ \dt
\eF
\[
- \int_0^\tau \intO{  \vr \Big( s (\vr, \vt) - s(\vrE, \vtE) \Big) \Big[ \partial_t \vtE
+ \vuE \cdot \Grad \vtE \Big] } \ \dt
\]
\[
+ \int_0^\tau \intO{ \left( \left( 1 - \frac{\vr}{\vrE} \right) \partial_t p(\vrE, \vtE) -
\frac{\vr}{\vrE} \vu \cdot \Grad p(\vrE, \vtE) \right) } \ \dt
\]
\[
\aleq \int_0^\tau \intO{ \ess{ p_M(\vrE, \vtE) - \frac{\partial p_M(\vrE, \vtE)}{\partial \vr} (\vrE - \vr)   - p_M(\vr, \vt) } \Div \vuE   } \ \dt
\]
\[
- \int_0^\tau \intO{  \vr \Big( s_M (\vr, \vt) - s_M (\vrE, \vtE) \Big) \ess{ \partial_t \vtE
+ \vuE \cdot \Grad \vtE } } \ \dt
\]
\[
+ \int_0^\tau \intO{(\vr - \vrE)  \vrE \frac{\partial s_M (\vrE, \vtE)}{\partial \vr} \ess{ \partial_t \vtE + \vuE \cdot \Grad \vtE } } \ \dt
\]
\[
+ \int_0^\tau \intO{ \mathcal{E} \left( \vr, \vt, \vu \Big| \vrE, \vtE, \vuE \right) } + \mathcal{O}(a)
\]
\[
\aleq \int_0^\tau \intO{ \ess{ p_M(\vrE, \vtE) - \frac{\partial p_M(\vrE, \vtE)}{\partial \vr} (\vrE - \vr)   - p_M(\vr, \vt) } \Div \vuE   } \ \dt
\]
\[
- \int_0^\tau \intO{  \vrE \Big( s_M (\vr, \vt) -  s_M (\vrE, \vtE) \Big) \ess{ \partial_t \vtE
+ \vuE \cdot \Grad \vtE } } \ \dt
\]
\[
+ \int_0^\tau \intO{(\vr - \vrE)  \vrE \frac{\partial s_M (\vrE, \vtE)} {\partial \vr} \ess{ \partial_t \vtE + \vuE \cdot \Grad \vtE } } \ \dt
\]
\[
+ \int_0^\tau \intO{ \mathcal{E} \left( \vr, \vt, \vu \Big| \vrE, \vtE, \vuE \right) } + \mathcal{O}(a)
\]
\[
\aleq \int_0^\tau \intO{ \ess{ p_M(\vrE, \vtE) - \frac{\partial p_M(\vrE, \vtE)}{\partial \vr} (\vrE - \vr)   - p_M(\vr, \vt) } \Div \vuE   } \ \dt
\]
\[
- \int_0^\tau \intO{  r \Big( s_M (\vr, \vt) -(\vr - \vrE) \frac{\partial s_M (\vrE, \vtE)}{\partial \vr}-  s_M (\vrE, \vtE) \Big) \ess{ \partial_t \vtE
+ \vuE \cdot \Grad \vtE } } \ \dt
\]
\[
+ \int_0^\tau \intO{ \mathcal{E} \left( \vr, \vt, \vu \Big| \vrE, \vtE, \vuE \right) } + \mathcal{O}(a)
\]
\[
\aleq \int_0^\tau \intO{ \ess{ p_M(\vrE, \vtE) - \frac{\partial p_M(\vrE, \vtE)}{\partial \vr} (\vrE - \vr)  - \frac{\partial p_M(\vrE, \vtE)}{\partial \vt} (\vtE - \vt) - p_M(\vr, \vt) } \Div \vuE   } \ \dt
\]
\[
- \int_0^\tau \int_\Omega  \vrE \Big( s_M (\vr, \vt) -(\vr - \vrE) \frac{\partial s_M (\vrE, \vtE)}{\partial \vr}
-(\vt - \vtE) \frac{\partial s_M (\vrE, \vtE)}{\partial \vt} -
s_M (\vrE, \vtE) \Big) \times
\]
\[
\times \ess{ \partial_t \vtE
+ \vuE \cdot \Grad \vtE }  \ \dx \ \dt
\]
\[
+ \int_0^\tau \intO{ \ess{ \frac{\partial p_M(\vrE, \vtE)}{\partial \vt} (\vtE - \vt) } \Div \vuE } \ \dt
\]
\[
- \int_0^\tau \intO{ \vrE (\vt - \vtE) \frac{\partial s_M (\vrE, \vtE)}{\partial \vt} \ess{ \partial_t \vtE
+ \vuE \cdot \Grad \vtE } } \ \dt
\]
\[
+ \int_0^\tau \intO{ \mathcal{E} \left( \vr, \vt, \vu \Big| \vrE, \vtE, \vuE \right) } + \mathcal{O}(a).
\]

Using the quadratic estimates (\ref{cc1}), the above simplifies to
\[
\int_0^\tau \intO{ \ess{ p_M(\vrE, \vtE) - \frac{\partial p_M(\vrE, \vtE)}{\partial \vr} (\vrE - \vr)  - \frac{\partial p_M(\vrE, \vtE)}{\partial \vt} (\vtE - \vt) - p_M(\vr, \vt) } \Div \vuE   } \ \dt
\]
\[
- \int_0^\tau \int_\Omega  \vrE \Big( s_M (\vr, \vt) -(\vr - \vrE) \frac{\partial s_M (\vrE, \vtE)}{\partial \vr}
-(\vt - \vtE) \frac{\partial s_M (\vrE, \vtE)}{\partial \vt} -
s_M (\vrE, \vtE) \Big) \times
\]
\[
\times \ess{ \partial_t \vtE
+ \vuE \cdot \Grad \vtE }  \ \dx \ \dt
\]
\[
+ \int_0^\tau \intO{ \ess{ \frac{\partial p_M(\vrE, \vtE)}{\partial \vt} (\vtE - \vt) } \Div \vuE } \ \dt
\]
\[
- \int_0^\tau \intO{ \vrE (\vt - \vtE) \frac{\partial s_M (\vrE, \vtE)}{\partial \vt} \ess{ \partial_t \vtE
+ \vuE \cdot \Grad \vtE } } \ \dt
\]
\[
+ \int_0^\tau \intO{ \mathcal{E} \left( \vr, \vt, \vu \Big| \vrE, \vtE, \vuE \right) } + \mathcal{O}(a)
\]
\[
\aleq
\int_0^\tau \intO{ \ess{ \frac{\partial p_M(\vrE, \vtE)}{\partial \vt} (\vtE - \vt) } \Div \vuE } \ \dt
\]
\[
+ \int_0^\tau \intO{ \ess{ (\vtE - \vt) \frac{\partial s_M (\vrE, \vtE)}{\partial \vt}  \left[ \partial_t (\vrE \vtE)
+ \Div (\vrE \vtE \vuE ) \right]} } \ \dt
\]
\[
+ \int_0^\tau \intO{ \mathcal{E} \left( \vr, \vt, \vu \Big| \vrE, \vtE, \vuE \right) } + \mathcal{O}(a).
\]

Finally, seeing that
\[
\vtE \frac{\partial s_M(\vrE, \vtE)}{\partial \vt} = \frac{\partial e_M(\vrE, \vtE)}{\partial \vt} = c_v (\vrE, \vtE),
\]
and, in accordance with (\ref{m3}),
\[
c_v (\vrE, \vtE) \left[ \partial_t (\vrE \vtE) + \Div (\vrE \vtE \vuE) \right] + \vtE \frac{\partial p_M(\vrE, \vtE)}{\partial \vt} \Div \vuE = 0,
\]
we may infer that
\bFormula{oi3}
\int_0^\tau \intO{ \vr \Big(  \partial_t \vuE +  \vuE \cdot \Grad \vuE \Big) \cdot (\vuE - \vu) } \ \dt
- \int_0^\tau \intO{ p(\vr, \vt) \Div \vuE   } \ \dt
\eF
\[
- \int_0^\tau \intO{  \vr \Big( s (\vr, \vt) - s(\vrE, \vtE) \Big) \Big[ \partial_t \vtE
+ \vuE \cdot \Grad \vtE \Big] } \ \dt
\]
\[
+ \int_0^\tau \intO{ \left( \left( 1 - \frac{\vr}{\vrE} \right) \partial_t p(\vrE, \vtE) -
\frac{\vr}{\vrE} \vu \cdot \Grad p(\vrE, \vtE) \right) } \ \dt
\]
\[
\aleq \int_0^\tau \mathcal{E} \left( \vr, \vt, \vu \Big| \vrE, \vtE, \vuE \right)  \ \dt + \mathcal{O}(a).
\]

Summing up the estimates (\ref{vhr3}), (\ref{oi3}) and applying Gronwall's lemma, we conclude
\[
\mathcal{E} \left( \vr, \vt, \vu \Big| \vrE, \vtE, \vuE \right)(\tau)
\]
\[
\aleq c(T,M,D) \left( \mathcal{E} \left( \vr_0, \vt_0, \vu_0 \Big| \vr_{0,E}, \vt_{0,E}, \vu_{0,E} \right)
+ \max \left\{a, \nu, \omega, \lambda,
\frac{\nu}{\sqrt{a}}, \frac{\omega}{a}, \left( \frac{a}{ \sqrt{\nu^3
\lambda} } \right)^{1/3} \right\} \right)
\]
for a.a. $\tau \in (0,T)$. This relation, together with the bounds established in (\ref{cc1}), (\ref{cc2}) give rise to
(\ref{m8}). We have proved Theorem \ref{Tm1} and Corollary \ref{Cm1}.

\section{Concluding remarks}
\label{F}

The piece of information provided by (\ref{konec}) is slightly better and definitely more explicit in terms of the convergence rate than stated
in Corollary \ref{Cm1}.
One could possibly do slightly better with respect to the values of the scaling parameters but it seems hard to remove the path dependence of the limit process.

The same result with obvious modifications in the proof can be obtained for a more general class of viscosity coefficients satisfying
\[
| \mu'(\vt) | \leq c, \ \underline{\mu} (1 + \vt^b) \leq \mu(\vt) \leq \Ov{\mu} (1 + \vt^b) ,\ 0 \leq \eta(\vt) \leq \Ov{\eta}(1 + \vt^b)
\  \mbox{for all}\ \vt \geq 0, \ \mbox{with}\ \frac{2}{5} < b \leq 1,
\]
in place of (\ref{H8}).

Similarly to the incompressible case, the situation described by the no-slip conditions
\[
\vu|_{\partial \Omega} = 0
\]
for the velocity in the Navier-Stokes-Fourier system remains an outstanding open problem.

\def\cprime{$'$} \def\ocirc#1{\ifmmode\setbox0=\hbox{$#1$}\dimen0=\ht0
  \advance\dimen0 by1pt\rlap{\hbox to\wd0{\hss\raise\dimen0
  \hbox{\hskip.2em$\scriptscriptstyle\circ$}\hss}}#1\else {\accent"17 #1}\fi}


\end{document}